\input amstex
\documentstyle{amsppt}
%
\catcode`@=11
\redefine\output@{%
  \def\break{\penalty-\@M}\let\par\endgraf
  \ifodd\pageno\global\hoffset=105pt\else\global\hoffset=8pt\fi  
  \shipout\vbox{%
    \ifplain@
      \let\makeheadline\relax \let\makefootline\relax
    \else
      \iffirstpage@ \global\firstpage@false
        \let\rightheadline\frheadline
        \let\leftheadline\flheadline
      \else
        \ifrunheads@ 
        \else \let\makeheadline\relax
        \fi
      \fi
    \fi
    \makeheadline \pagebody \makefootline}%
  \advancepageno \ifnum\outputpenalty>-\@MM\else\dosupereject\fi
}
\def\Beta{\mathchar"0\hexnumber@\rmfam 42}
\catcode`\@=\active
\nopagenumbers
\chardef\textvolna='176

\chardef\bigalpha='013

\accentedsymbol\DD{D\kern -1.5pt D}

\chardef\degree="5E
\def\compos{\,\raise 1pt\hbox{$\sssize\circ$} \,}

\def\blue#1{#1}

\catcode`#=11\def\diez{#}\catcode`#=6
\catcode`&=11\catcode`&=4
\catcode`_=11\def\podcherkivanie{_}\catcode`_=8
\catcode`\^=11\catcode`\^=7
\catcode`~=11\catcode`~=\active
\def\mycite#1{\cite{\blue{#1}}\immediate\special{ps:
     ShrHPSdict begin /ShrBORDERthickness 0 def}}

\def\mytag#1{%
    \tag#1}
\def\mythetag#1{\thetag{\blue{#1}}\immediate\special{ps:
     ShrHPSdict begin /ShrBORDERthickness 0 def}}
\def\myrefno#1{\no#1}
\def\myhref#1#2{\blue{#2}\immediate\special{ps:
     ShrHPSdict begin /ShrBORDERthickness 0 def}}

\def\mytheorem#1{\csname proclaim\endcsname{Theorem #1}}
\def\mytheoremwithtitle#1#2{\csname proclaim\endcsname{Theorem #1#2}}
\def\mythetheorem#1{\blue{#1}\immediate\special{ps:
     ShrHPSdict begin /ShrBORDERthickness 0 def}}
\def\mylemma#1{\csname proclaim\endcsname{Lemma #1}}
\def\mylemmawithtitle#1#2{\csname proclaim\endcsname{Lemma #1#2}}

\def\mycorollary#1{\csname proclaim\endcsname{Corollary #1}}

\def\mydefinition#1{\definition{Definition #1}}
\def\mythedefinition#1{\blue{#1}\immediate\special{ps:
     ShrHPSdict begin /ShrBORDERthickness 0 def}}
\def\myconjecture#1{\csname proclaim\endcsname{Conjecture #1}}
\def\myconjecturewithtitle#1#2{\csname proclaim\endcsname{Conjecture #1#2}}

\def\myproblem#1{\csname proclaim\endcsname{Problem #1}}
\def\myproblemwithtitle#1#2{\csname proclaim\endcsname{Problem #1#2}}


\pagewidth{360pt}
\pageheight{606pt}
\topmatter
\title
Multiple discriminants and critical values of a multivariate polynomial.
\endtitle
\rightheadtext{Multiple discriminants and critical values \dots}
\author
Ruslan Sharipov
\endauthor
\address Bashkir State University, 32 Zaki Validi street, 450074 Ufa, Russia
\endaddress
\email
\myhref{mailto:r-sharipov\@mail.ru}{r-sharipov\@mail.ru}
\endemail
\abstract
     A critical value of a function is the value of this function at one of its 
critical points. Each critical point of a differentiable multivariate function 
is described by the equations which consist in equating to zero all of its 
partial derivatives. However, in general case there is no equation for the 
corresponding critical value. The case of polynomials is different. In the present 
paper an equation for critical values of a polynomial is derived. 
\endabstract
\endtopmatter
\TagsOnRight
\document

\head
1. Introduction.
\endhead
     Let $f(x_1,\ldots,x_n)$ be a smooth real valued multivariate function
in $\Bbb R^n$ or in some open domain $\Omega\subset\Bbb R^n$. Critical points 
of  the function $f$ are determined by solving the following system of equations
with respect to the variables $x_1,\,\ldots,\,x_n$:
$$
\hskip -2em
\cases\dfrac{\partial f(x_1,\ldots,x_n)}{\partial x_1}=0,\\
\ .\ .\ .\ .\ .\ .\ .\ .\ .\ .\ .\ .\ .\ .\\
\dfrac{\partial f(x_1,\ldots,x_n)}{\partial x_n}=0.
\endcases
\mytag{1.1}
$$
They can be maxima, minima, saddle points, and other types of critical points. 
\par
     Let $v$ be the value of $f(x_1,\ldots,x_n)$ at one of its critical points
given by the equations \mythetag{1.1}. Then we can extend the system of equations
\mythetag{1.1} as follows:
$$
\hskip -2em
\cases\dfrac{\partial f(x_1,\ldots,x_n)}{\partial x_1}=0,\\
\ .\ .\ .\ .\ .\ .\ .\ .\ .\ .\ .\ .\ .\ .\\
\dfrac{\partial f(x_1,\ldots,x_n)}{\partial x_n}=0,\\
\vspace{1ex}
f(x_1,\ldots,x_n)-v=0.
\endcases
\mytag{1.2}
$$
The parameter $v$ in \mythetag{1.2} is treated as a new variable independent of $x_1\,\ldots,\,x_n$, i\.\,e\. \mythetag{1.2} is a system of $n+1$ equations for
$n+1$ variables. Theoretically, one can eliminate the variables $x_1\,\ldots,\,x_n$
from the system \mythetag{1.2} thus producing one equation for one variable $v$,
which is the critical value of $f$:
$$
\pagebreak 
\hskip -2em
F(v)=0.
\mytag{1.3}
$$
In practice it is rather difficult to derive such an equation. In the present 
paper we derive the equation of the form \mythetag{1.3} for the case where 
$f(x_1,\ldots,x_n)$ is a multivariate polynomial. For the sake of simplicity,
from now on we treat $x_1\,\ldots,\,x_n$ as complex variables so that they
constitute a point in $\Bbb C^{\kern 0.5pt n}$.\par
\head
2. The case of a univariate polynomial. 
\endhead
     Let's begin with the case where $n=1$. Denoting $x_1=x$ for the sake 
of simplicity, assume that $f(x)$ is a univariate polynomial of $n$-th degree:
$$
f(x)=a_0\,x^n+a_1\,x^{n-1}+\ldots+a_{n-1}\,x+a_n.
$$
Then $p(x)=f(x)-v$ is also a univariate polynomial of $n$-th degree: 
$$
\hskip -2em
p(x)=a_0\,x^n+a_1\,x^{n-1}+\ldots+a_{n-1}\,x+(a_n-v).
\mytag{2.1}
$$
Using $p(x)$, in this case we can write the system of equations \mythetag{1.2} 
as follows:
$$
\cases p^{\kern 1pt\prime}(x)=0,\\
p(x)=0.
\endcases
$$
It is known that a univariate polynomial $p(x)$ vanishes along with its first
derivative $p^{\kern 1pt\prime}(x)$ at some point $x\in\Bbb C$ if and only its
discriminant is zero (see \mycite{1}):
$$
\hskip -2em
D_p=0.
\mytag{2.2}
$$
The coefficients of the polynomial $p(x)$ in \mythetag{2.1} depend on $v$. 
Therefore the discriminant $D_p$ in \mythetag{2.2} depends on $v$. And we
write \mythetag{2.2} as 
$$
\hskip -2em
D_p(v)=0.
\mytag{2.3}
$$
The equation \mythetag{2.3} is a required equation of the form \mythetag{1.3}. 
\mytheorem{2.1} A complex number $v$ is a critical value of a univariate  polynomial 
$f(x)$ if and only if it is a root of the equation \mythetag{2.3}, where $D_p(v)$
is the discriminant of the polynomial $p(x)=f(x)-v$. 
\endproclaim
     Theorem~2.1 follows immediately from the above mentioned basic property of
the discriminant of a univariate polynomial.\par
\head
3. The case of a multivariate polynomial. 
\endhead
     This case is somewhat similar to the previous one. Here we can replace
the multivariate polynomial $f(x_1,\ldots,x_n)$ with the polynomial 
$$
\hskip -2em
p(x_1,\ldots,x_n)=f(x_1,\ldots,x_n)-v. 
\mytag{3.1}
$$
Due to \mythetag{3.1} we can write the equations \mythetag{1.2} as follows:
$$
\pagebreak 
\hskip -2em
\cases\dfrac{\partial p(x_1,\ldots,x_n)}{\partial x_1}=0,\\
\ .\ .\ .\ .\ .\ .\ .\ .\ .\ .\ .\ .\ .\ .\\
\dfrac{\partial p(x_1,\ldots,x_n)}{\partial x_n}=0,\\
\vspace{1ex}
p(x_1,\ldots,x_n)=0.
\endcases
\mytag{3.2}
$$
The concept of the discriminant of a multivariate polynomial is not 
commonly known. Its definition can be found in \mycite{2} (see also
\mycite{3} and \mycite{4}).	
\mydefinition{3.1} The discriminant $D_p$ of a multivariate polynomial 
$p(x_1,\ldots,x_n)$ is a polynomial function of the coefficients of $p$ 
such that its own coefficients are integer, such that it is irreducible 
over $\Bbb Z$, and such that $D_p=0$ if and only if the system of 
equations \mythetag{3.2} has an least one solution in complex 
numbers $x_1,\,\ldots,\,x_n$. 
\enddefinition
     As it was said in \mycite{2} with the reference to \mycite{5},
discriminants of multivariate polynomials were first considered 
by G.~Boole. The reference \mycite{5} is taken from \mycite{2} ``as is''.
It looks quite uncertain, since no information on publishers is 
provided (see more detailed historical research in \mycite{6} and 
\mycite{7}).\par
     The coefficients of the polynomial $p(x_1,\ldots,x_n)$ in \mythetag{3.1} 
depend on $v$. Therefore the discriminant $D_p$ depends on $v$. Hence we can
write the equation
$$
\hskip -2em	
D_p(v)=0.
\mytag{3.3}
$$
The equation \mythetag{3.3} is of the form \mythetag{1.3}. It solves our 
problem of reducing \mythetag{1.2} to a single equation for $v$ through the
following theorem which is immediate from the above 
Definition~\mythedefinition{3.1}.\par
\mytheorem{3.1} A complex number $v$ is a critical value of a multivariate  
polynomial $f(x_1,\ldots,x_n)$ if and only if it is a root of the equation 
\mythetag{3.3}, where $D_p(v)$ is the discriminant of the polynomial 
\mythetag{3.1}.
\endproclaim
     However, the matter is that there is no simple formula for the 
discriminant $D_p$ in the case of a multivariate polynomial $p$. For this
reason below we consider a different approach to deriving an equation
of the form \mythetag{1.3} from \mythetag{1.2}.\par 
\head
4. Critical points of discriminants. 
\endhead
    Let $p(x,y)$ be a univariate polynomial of $n$-th degree in $x$, i\.\,e\.
$$
\hskip -2em
p(x,y)=a_n(y)\,x^n+a_{n-1}(y)\,x^{n-1}+\ldots+a_1(y)\,x+a_0(y),
\mytag{4.1}
$$
whose coefficients are smooth function of a complex variable $y\in\Bbb C$.
Assume that at some point $(x_0,y_0)\in\Bbb C^{\kern 0.5pt 2}$, where 
$a_n(y_0)\neq 0$, the following equations are fulfilled:
$$
\xalignat 3
&\hskip -2em
p^{\kern 1pt\prime}_x(x_0,y_0)=0,
&&p^{\kern 1pt\prime}_y(x_0,y_0)=0,
&&p(x_0,y_0)=0.
\mytag{4.2}
\endxalignat 
$$
\mytheorem{4.1} If $p$ is a univariate polynomial of the form 
\mythetag{4.1} whose coefficients are smooth function of a complex
variable $y$ and if the equations \mythetag{4.2} are fulfilled at
some point $(x_0,y_0)\in\Bbb C^{\kern 0.5pt 2}$, where $a_n(y_0)\neq 0$, 
then they imply the equations 
$$
\xalignat 2
&\hskip -2em
D^{\kern 1pt\prime}(y_0)=0,
&&D(y_0)=0
\mytag{4.3}
\endxalignat 
$$
for the discriminant $D$ of $p$ which are fulfilled at the point 
$y_0\in\Bbb C$, where $y_0$ is the second coordinate of the corresponding 
point $(x_0,y_0)\in\Bbb C^{\kern 0.5pt 2}$. 
\endproclaim
     Since $a_n(y_0)\neq 0$ in Theorem~\mythetheorem{4.1}, we can divide 
the polynomial \mythetag{4.1} by $a_n(y)$ and proceed to the following 
monic polynomial:
$$
\hskip -2em
\tilde p(x,y)=x^n+\frac{a_{n-1}(y)}{a_n(y)}\,x^{n-1}+\ldots
+\frac{a_1(y)}{a_n(y)}\,x
+\frac{a_0(y)}{a_n(y)}.
\mytag{4.4}
$$
The discriminants of the polynomials \mythetag{4.1} and \mythetag{4.4} are 
related as follows:
$$
\hskip -2em
D_p(y)=(a_n(y))^{2\,n-2}\,D_{\tilde p}(y).
\mytag{4.5}
$$
Relying on \mythetag{4.5}, we can write the monic polynomial 
$$
\hskip -2em
p(x,y)=x^n+a_{n-1}(y)\,x^{n-1}+\ldots+a_1(y)\,x+a_0(y)
\mytag{4.6}
$$
and then reformulate Theorem~\mythetheorem{4.1} in the following way. 
\mytheorem{4.2} If $p$ is a monic univariate polynomial of the form 
\mythetag{4.6} whose coefficients are smooth function of a complex
variable $y$ and if the equations \mythetag{4.2} are fulfilled at
some point $(x_0,y_0)\in\Bbb C^{\kern 0.5pt 2}$, then they imply the 
equations \mythetag{4.3} for the discriminant $D$ of $p$ which are 
fulfilled at the point $y_0\in\Bbb C$, where $y_0$ is the second 
coordinate of the corresponding point $(x_0,y_0)\in\Bbb C^{\kern 0.5pt 2}$. 
\endproclaim
    Theorems~\mythetheorem{4.1} and \mythetheorem{4.2} are equivalent
to each other due to the relationship \mythetag{4.5}.\par
\head
5. Proof of Theorem~\mythetheorem{4.2} in the case of a quadratic
polynomial. 
\endhead
     Assume that $n=2$ in \mythetag{4.6}. Then $p(x,y)$ is a quadratic 
polynomial:
$$
\hskip -2em
p(x,y)=x^2+a_1(y)\,x+a_0(y).
\mytag{5.1}
$$
One can easily calculate the partial derivatives of $p(x,y)$:
$$
\xalignat 2
&\hskip -2em
p^{\kern 1pt\prime}_x(x,y)=2\,x+a_1(y),
&&p^{\kern 1pt\prime}_y(x,y)=a^{\kern 1pt\prime}_1(y)\,x+
a^{\kern 1pt\prime}_0(y). 
\mytag{5.2}
\endxalignat
$$
Substituting \mythetag{5.1} and \mythetag{5.2} into \mythetag{4.2},
we derive a system of three equations:
$$
\hskip -2em
\cases
2\,x_0+a_1(y_0)=0,\\
a^{\kern 1pt\prime}_1(y_0)\,x_0+a^{\kern 1pt\prime}_0(y_0)=0,\\
(x_0)^2+a_1(y_0)\,x_0+a_0(y_0)=0.
\endcases
\mytag{5.3}
$$
The discriminant of the polynomial \mythetag{5.1} is calculated as follows:
$$
\hskip -2em
D(y)=(a_1(y))^2-4\,a_0(y). 
\mytag{5.4}
$$
Its derivative $D^{\kern 1pt\prime}(y)$ is also easily calculated:
$$
D^{\kern 1pt\prime}(y)=2\,a_1(y)\,a^{\kern 1pt\prime}_1(y)
-4\,a^{\kern 1pt\prime}_0(y). 
\mytag{5.5}
$$\par
     Note that the first equation in \mythetag{5.3} can be resolved with 
respect to $x_0$:
$$
\hskip -2em
x_0=-\frac{a_1(y_0)}{2}.
\mytag{5.6}
$$
Substituting \mythetag{5.6} into the third equation \mythetag{5.3} we
derive  
$$
\pagebreak
\hskip -2em
-\frac{(a_1(y_0))^2}{4}+a_0(y_0)=0.
\mytag{5.7}
$$
Comparing \mythetag{5.7} with \mythetag{5.4}, we see that \mythetag{5.3}
implies 
$$
\hskip -2em
D(y_0)=0.
\mytag{5.8}
$$
Then we substitute \mythetag{5.6} into the second equation \mythetag{5.3}.
As a result we get
$$
\hskip -2em
-\frac{\,a_1(y_0)\,a^{\kern 1pt\prime}_1(y_0)}{2}
+\,a^{\kern 1pt\prime}_0(y_0)=0.
\mytag{5.9}
$$
Comparing \mythetag{5.9} with \mythetag{5.5}, we see that \mythetag{5.3}
implies 
$$
\hskip -2em
D^{\kern 1pt\prime}(y_0)=0.
\mytag{5.10}
$$
The rest is to note that \mythetag{5.10} and \mythetag{5.8} coincide
with \mythetag{4.3} and conclude that \mythetag{4.2} implies \mythetag{4.3}
in the case of the polynomial \mythetag{5.1}. Thus, Theorem~\mythetheorem{4.2}
is proved for the case of quadratic polynomials. 
\head
6. Proof of Theorem~\mythetheorem{4.2} in the case of a 
double root polynomial.
\endhead
     If $y=y_0$ is fixed, the equations \mythetag{4.2} mean that the polynomial 
$p$ vanishes along with its derivative $p^{\kern 1pt\prime}_x$ at the point
$x=x_0$, i\.\,e\. it has a multiple root at this point. The case of a double
root is the most simple in this situation.\par
     Let $x_1\,\ldots,\,x_n$ be roots of the polynomial \mythetag{4.6}. They
depend on $y$, i\.\,e\.
$$
\hskip -2em
x_i=x_i(y),\ \ i=1,\,\ldots,\,n.
\mytag{6.1}
$$
It is known that roots of a univariate polynomial are continuous functions
of its coefficients (see \mycite{8}). Hence in our case the roots \mythetag{6.1}
are continuous functions of $y$. Two of them tend to $x_0$ as $y\to y_0$. 
Without loss of generality we can set
$$
\hskip -2em
x_1(y)\to x_0\text{\ \ and \ }x_2(y)\to x_0\text{\ \ as  \ }
y\to y_0. 
\mytag{6.2}
$$
Assume that the roots $x_3(y_0),\,\ldots,\,x_n(y_0)$ are distinct and different
from the double root $x_1(y_0)=x_2(y_0)=x_0$. Under this assumption the roots
$$
\hskip -2em
x_3(y),\,\ldots,\,x_n(y)
\mytag{6.3}
$$
are smooth functions of $y$ in some neighborhood of the point $y=y_0$ 
(see \mycite{8}). Unlike them the roots $x_1(y)$ and $x_2(y)$ in \mythetag{6.2} 
are not necessarily smooth, though they are continuous functions of $y$ in this 
neighborhood of the point $y=y_0$.\par
    Using the roots \mythetag{6.2} and \mythetag{6.3}, we define two 
polynomials which are two complementary factors of the initial polynomial 
$p(x,y)$ in \mythetag{4.6}:
$$
\xalignat 2
&\hskip -2em
q(x,y)=\prod_{i=1}^2(x-x_i(y)),
&&r(x,y)=\prod_{i=3}^n(x-x_i(y)).
\mytag{6.4}
\endxalignat
$$
Indeed, from \mythetag{6.4} we easily derive the equality 
$$
\hskip -2em
p(x,y)=q(x,y)\,r(x,y).
\mytag{6.5}
$$\par
    The coefficients of the polynomial $r(x,y)$ in \mythetag{6.4} are
smooth functions of $y$ since they are elementary symmetric functions of
the smooth roots \mythetag{6.3} (see \mycite{1}). The polynomial $q(x,y)$
in \mythetag{6.4} is quadratic: 
$$
\hskip -2em
q(x,y)=x^2+\beta_1(y)\,x+\beta_0(y). 
\mytag{6.6}
$$
Its coefficients are also smooth functions of $y$ since due to \mythetag{6.5}
the polynomial \mythetag{6.6} can be produced by means of the polynomial 
long division algorithm (see \mycite{9}) with the monic polynomial
$r(x,y)$ as a divisor:
$$
q(x,y)=p(x,y)\div r(x,y).
$$\par
     Let's recall that the roots $x_3(y_0),\,\ldots,\,x_n(y_0)$ are different
from the double root $x_0$. Therefore from \mythetag{6.4} we derive the
inequality
$$
\hskip -2em
r(x_0,y_0)=\prod_{i=3}^n(x_0-x_i(y_0))\neq 0. 
\mytag{6.7}
$$
Similarly, applying \mythetag{6.2} to $q(x,y)$ in \mythetag{6.4}, we derive 
the equality
$$
\hskip -2em
q(x_0,y_0)=0.
\mytag{6.8}
$$
Differentiating the equality \mythetag{6.5}, we find that 
$$
\hskip -2em
\aligned
&p^{\kern 1pt\prime}_x(x_0,y_0)=
q^{\kern 1pt\prime}_x(x_0,y_0)\,r(x_0,y_0)+
q(x_0,y_0)\,r^{\kern 1pt\prime}_x(x_0,y_0),\\
&p^{\kern 1pt\prime}_y(x_0,y_0)=
q^{\kern 1pt\prime}_y(x_0,y_0)\,r(x_0,y_0)+
q(x_0,y_0)\,r^{\kern 1pt\prime}_y(x_0,y_0).
\endaligned
\mytag{6.9}
$$
Now, applying \mythetag{6.8} to \mythetag{6.9}, we obtain the
equalities
$$
\hskip -2em
\aligned
&p^{\kern 1pt\prime}_x(x_0,y_0)=q^{\kern 1pt\prime}_x(x_0,y_0)
\,r(x_0,y_0),\\
&p^{\kern 1pt\prime}_y(x_0,y_0)=q^{\kern 1pt\prime}_y(x_0,y_0)
\,r(x_0,y_0).
\endaligned
\mytag{6.10}
$$
If the equations \mythetag{4.2} are fulfilled, then,
using \mythetag{6.7}, from \mythetag{6.10} we derive
$$
\xalignat 2
&\hskip -2em
q^{\kern 1pt\prime}_x(x_0,y_0)=0,
&&q^{\kern 1pt\prime}_y(x_0,y_0)=0.
\mytag{6.11}
\endxalignat 
$$\par
     The equalities \mythetag{6.11} combined with \mythetag{6.8}
mean that the equations \mythetag{4.2} for $p$, once they are fulfilled, 
imply similar equations 
$$
\xalignat 3
&\hskip -2em
q^{\kern 1pt\prime}_x(x_0,y_0)=0,
&&q^{\kern 1pt\prime}_y(x_0,y_0)=0,
&&q(x_0,y_0)=0
\mytag{6.12}
\endxalignat 
$$
for the quadratic polynomial \mythetag{6.6} whose coefficients are
smooth functions of $y$. Theorem~\mythetheorem{4.2} is already proved
for quadratic polynomials. Therefore, applying this theorem to $q(x,y)$
and \mythetag{6.12}, we derive the equalities 
$$
\xalignat 2
&\hskip -2em
D^{\kern 1pt\prime}_q(y_0)=0,
&&D_q(y_0)=0,
\mytag{6.13}
\endxalignat 
$$
where $D_q(y)$ is the discriminant of the polynomial $q(x,y)$
in \mythetag{6.6}.\par
     Let's proceed to the discriminant $D(y)$ of the polynomial $p(x,y)$
in \mythetag{4.6}. In terms of its roots the discriminant $D(y)$ is given 
by the formula
$$
\hskip -2em
D(y)=\prod^n_{i<j}(x_i(y)-x_j(y))^2
\mytag{6.14}
$$
(see \mycite{1}). If we denote through $D_r(y)$ the discriminant of the
polynomial $r(x,y)$ in \mythetag{6.4}, then we can factorize the discriminant
\mythetag{6.14} as follows:
$$
\hskip -2em
D(y)=D_q(y)\,D_r(y)\,\prod^2_{i=1}\prod^n_{j=3}(x_i(y)-x_j(y))^2.
\mytag{6.15}
$$
Taking into account the formula for $q(x,y)$ in \mythetag{6.4}, we can transform
\mythetag{6.15} as 
$$
\hskip -2em
D(y)=D_q(y)\,D_r(y)\biggl(\,\prod^n_{j=3}q(x_j(y),y)\biggr)^{\kern -1pt 2}.
\mytag{6.16}
$$\par
     Relying on \mythetag{6.16}, lets introduce the following notation:
$$
\hskip -2em
\alpha(y)=D_r(y)\,\biggl(\,\prod^n_{j=3}q(x_j(y),y)\biggr)^{\kern -1pt 2}.
\mytag{6.17}
$$
The function $\alpha(y)$ in \mythetag{6.17} is a smooth function of $y$ since
the coefficients of the polynomial $q(x,y)$ in \mythetag{6.6} and the roots
$x_3(y),\,\ldots,\,x_n(y)$ of the polynomial $r(x,y)$ in \mythetag{6.3} 
are smooth functions of $y$. Applying \mythetag{6.17} to \mythetag{6.16}, we
get
$$
\hskip -2em
D(y)=D_q(y)\,\alpha(y).
\mytag{6.18}
$$
Both multiplicands $D_q(y)$ and $\alpha(y)$ in \mythetag{6.18} are smooth
functions of $y$. Therefore, differentiating \mythetag{6.18}, we derive 
the following formula:
$$
\hskip -2em
D^{\kern 1pt\prime}(y_0)=D^{\kern 1pt\prime}_q(y_0)\,\alpha(y_0)
+D_q(y_0)\,\alpha^{\kern 1pt\prime}(y_0).
\mytag{6.19}
$$ 
The rest is to apply \mythetag{6.13} to \mythetag{6.18} and \mythetag{6.19}
and derive \mythetag{4.3}. Thus we have proved that \mythetag{4.2} implies
\mythetag{4.3} in our present case. Theorem~\mythetheorem{4.2} is proved
in the case of a polynomial \mythetag{4.6} with exactly one double root
and the other simple roots. 
\head
7. Proof of Theorem~\mythetheorem{4.2} in general case.
\endhead
    Note that the polynomial $p(x,y)$ in \mythetag{4.6} depends on $y$ through
its coefficients. Therefore the derivative $D^{\kern 1pt\prime}(y)$ of its 
discriminant is calculated is follows:
$$
\hskip -2em
D^{\kern 1pt\prime}(y)=\sum^{n-1}_{i=0}\frac{\partial D}{\partial a_i}
\ a^{\kern 1pt\prime}_i(y).
\mytag{7.1}
$$
Let's denote $\delta p(x,y)=p^{\kern 1pt\prime}_y(x,y)$. Then, differentiating
\mythetag{4.6}, we get
$$
\hskip -2em
\delta p(x,y)=a^{\kern 1pt\prime}_{n-1}(y)\,x^{n-1}+\ldots+a^{\kern 1pt\prime}_1(y)
\,x+a^{\kern 1pt\prime}_0(y). 
\mytag{7.2}
$$
Let's denote $a_i=a^{\kern 1pt\prime}_i(y_0)$, $b_i=a^{\kern 1pt\prime}_i(y_0)$ 
and $\delta D=D^{\kern 1pt\prime}(y_0)$. Then, substituting $y=y_0$ into 
\mythetag{4.6}, \mythetag{7.2} and \mythetag{7.1}, we obtain two polynomials
$$
\aligned
&p(x)=x^n+a_{n-1}\,x^{n-1}+\ldots+a_1\,x+a_0,\\
\delta &p(x)=b_{n-1}\,x^{n-1}+\ldots+b_1\,x+b_0 
\endaligned
\mytag{7.3}
$$
with purely numeric coefficients and the purely numeric quantity 
$$
\hskip -2em
\delta D=\sum^{n-1}_{i=0}\frac{\partial D}{\partial a_i}\ b_i
\mytag{7.4}
$$
associated with the polynomials \mythetag{7.3}. The following definition is
terminological. It is designed for the sake of beauty. 
\mydefinition{7.1} The polynomial $\delta p(x)$ in \mythetag{7.3} is called 
a first variation of the polynomial $p(x)$, while the numeric quantity 
\mythetag{7.4} is called the first variation of the discriminant $D$ of $p$
associated with $\delta p(x)$. 
\enddefinition
    It turns out that the functional nature of the coefficients of $p(x,y)$
in \mythetag{4.6} is inessential in Theorem~\mythetheorem{4.2}. This theorem
can be reformulated as follows. 
\mytheorem{7.1} If a monic polynomial $p(x)$ in \mythetag{7.3} has a multiple
root of multiplicity $m\geqslant 2$ and if it shares this root with its 
first variation $\delta p(x)$, then the discriminant $D$ of $p(x)$ vanishes
along with its first variation $\delta D$ in \mythetag{7.4}. 
\endproclaim 
     Theorem~\mythetheorem{7.1} is equivalent to Theorem~\mythetheorem{4.2}.
The equivalence can be established using the linear functions $a_i(y)=a_i
+b_i\,(y-y_0)$ in \mythetag{4.6}. The result of the previous section means that 
we have already proved Theorem~\mythetheorem{7.1} for any monic polynomial $p(x)$ 
with exactly one double root and the other simple roots.\par
      Let $p(x)$ be a monic polynomial with at least one 
multiple root $x_0$. Then it can be produced as a limit of some sequence of monic 
polynomials $p_s(x)$ with exactly one double root $x=x_0$ and the other simple
roots. Assume that $p(x)$ shares its multiple root $x=x_0$ with its first 
variation $\delta p(x)$ in \mythetag{7.3}. Then
$$
\hskip -2em
\delta p(x_0)=0.
\mytag{7.5}
$$
The equality \mythetag{7.5} is written as a linear relationships with respect to 
the coefficients $b_0,\,\ldots,\,b_{n-1}$ of the polynomial $\delta p(x)$ in 
\mythetag{7.3}: 
$$
\hskip -2em
x^{n-1}_0\,b_{n-1}+\ldots+x_0\,b_1+b_0=0.
\mytag{7.6}
$$
The coefficients of the linear combination \mythetag{7.6} depend only on the root 
$x_0$, which is common for $p(x)$ and for any polynomial in the sequence $p_s(x)$.
This means that $p(x)$ and the polynomials $p_s(x)$ share the same first variation
$\delta p(x)$ obeying the relationship \mythetag{7.5}. This first variation
$\delta p(x)$ combined with each $p_s(x)$ according to \mythetag{7.4} produces a 
numeric sequence $\delta D_s$. The discriminant $D$ and its partial derivatives
in \mythetag{7.4} are smooth functions of $a_0,\,\ldots,\,a_{n-1}$. Therefore we 
have
$$
\hskip -2em
\delta D=\lim_{s\to\infty}\delta D_s. 
\mytag{7.7}
$$
The rest is to apply Theorem~\mythetheorem{7.1} to each $p_s(x)$ combined with
$\delta p(x)$. \pagebreak This yields $\delta D_s=0$. Substituting $\delta D_s=0$ 
into \mythetag{7.7}, we derive the required result $\delta D=0$. Thus, 
Theorem~\mythetheorem{7.1} is proved in its full generality. The same is true 
for Theorems~\mythetheorem{4.2} and \mythetheorem{4.1}, which are equivalent 
to Theorem~\mythetheorem{7.1}. 
\head
8. Application to critical values.
\endhead
    Having been equipped with Theorem~\mythetheorem{4.1}, now we return 
to our initial problem of deriving an equation of the form \mythetag{1.3}
from the equations \mythetag{1.2} in the case of a multivariate 
polynomial $f(x_1,\ldots,x_n)$. Passing from $f(x_1,\ldots,x_n)$ to the 
polynomial $p(x_1,\ldots,x_n)$ in \mythetag{3.1}, we treat $p(x_1,\ldots,x_n)$
as a univariate polynomial with respect to the variable $x=x_n$ and treat 
the other variables $x_1,\,\ldots,\,x_{n-1}$ and $v$ as parameters. Then we 
calculate the univariate discriminant $D$ of the polynomial $p(x_1,\ldots,x_n)$ 
with respect to the last variable $x=x_n$:
$$
\tilde p(x_1,\ldots,x_{n-1})=D(x_n,p(x_1,\ldots,x_n)).
\mytag{8.1}
$$
The univariate discriminant $D$ in \mythetag{8.1} acts as a nonlinear operator
sending the $n$-variate polynomial $p(x_1,\ldots,x_n)$ to the 
$(n-1)$\kern 1pt-variate polynomial $\tilde p(x_1,\ldots,x_{n-1})$. Setting 
$y=x_i$ and applying Theorem~\mythetheorem{4.1} for each $i=1,\,\ldots,\,n-1$, 
we derive the following system of equations from the equations \mythetag{3.2}:
$$
\hskip -2em
\cases\dfrac{\partial\tilde p(x_1,\ldots,x_{n-1})}{\partial x_1}=0,\\
\ .\ .\ .\ .\ .\ .\ .\ .\ .\ .\ .\ .\ .\ .\\
\dfrac{\partial\tilde p(x_1,\ldots,x_{n-1})}{\partial x_{n-1}}=0,\\
\vspace{1ex}
\tilde p(x_1,\ldots,x_{n-1})=0.
\endcases
\mytag{8.2}
$$
The structure of the equations \mythetag{8.2} is the same as the structure of
the equations \mythetag{3.2}. Therefore we can apply the operator \mythetag{8.1}
repeatedly:
$$
\hskip -2em
\DD_p=D(x_1,D(x_2,\ldots, D(x_n,p(x_1,\ldots,x_n))\ldots)).
\mytag{8.3}
$$ 
\mydefinition{8.1} The numeric quantity $\DD_p$ introduced through the formula
\mythetag{8.3} is called the multiple discriminant of a multivariate polynomial
$p(x_1,\ldots,x_n)$. 
\enddefinition
    Generally speaking, the multiple discriminant $\DD_p$ is different from the
discriminant $D_p$ of a multivariate polynomial introduced in 
Definition~\mythedefinition{3.1}. Unlike $D_p$, the formula \mythetag{8.3} 
provides a clear algorithm for calculating $\DD_p$.\par
     The coefficients of the polynomial $p(x_1,\ldots,x_n)$ in \mythetag{3.1} 
depend on $v$. Therefore the multiple discriminant $\DD_p$ depends on $v$. Hence 
we can write the equation
$$
\hskip -2em	
\DD_p(v)=0,
\mytag{8.4}
$$
where $\DD_p(v)=D(x_1,D(x_2,\ldots, D(x_n,f(x_1,\ldots,x_n)-v))\ldots))$.
\mytheorem{8.1} If a complex number $v$ is a critical value of a multivariate  
polynomial $f(x_1,\ldots,x_n)$, then it is a root of the equation \mythetag{8.4}.
\endproclaim
     Theorem~\mythetheorem{8.1} is proved by applying 
Theorem~\mythetheorem{4.1} repeatedly. Theorem~\mythetheorem{8.1} is similar to Theorem~\mythetheorem{3.1}, \pagebreak but it is somewhat weaker, i\.\,e\. each 
critical value of a multivariate polynomial $f(x_1,\ldots,x_n)$ is a root of the 
equation \mythetag{8.4}, but not each root of the equation \mythetag{8.4} is 
a critical value of $f(x_1,\ldots,x_n)$.\par 
\head
9. Conclusions.
\endhead
     Theorem~\mythetheorem{8.1} along with the equation \mythetag{8.4} and 
the formula \mythetag{8.3} constitutes the main result of the present paper. 
It can be applied in testing positivity of multivariate quartic forms and
forms of higher degrees. 

\enddocument
\end